# J. L. DOOB:
# FOUNDATIONS OF STOCHASTIC PROCESSES
# AND PROBABILISTIC POTENTIAL THEORY

By Ronald Getoor

*University of California at San Diego*

During the three decades from 1930 to 1960 J. L. Doob was, with the possible exception of Kolmogorov, the man most responsible for the transformation of the study of probability to a mathematical discipline. His accomplishments were recognized by both probabilists and other mathematicians in that he is the only person ever elected to serve as president of both the IMS and the AMS. This article is an attempt to discuss his contributions to two areas in which his work was seminal, namely, the foundations of continuous parameter stochastic processes and probabilistic potential theory.

**1. Introduction.** J. L. Doob was a true mathematical pioneer. He was instrumental in establishing and developing a number of areas of *mathematics* that became major topics for further research in the mathematical theory of probability. I have emphasized mathematics in the preceding sentence because Doob himself insisted on it. In fact, I think that what we owe most to him is his steadfast, unrelenting insistence that probability is mathematics. He clearly stated this in the preface to his 1953 book [36]—such numbers refer to Doob's publications in his accompanying complete bibliography that was compiled by D. Burkholder—where he wrote, "Probability is simply a branch of measure theory, with its own special emphasis and field of application ...." In one of his earliest papers on probability [10] he wrote, and the emphasis is his, "It is the purpose of this paper to show that *the statement that a successful system is impossible in a game of chance corresponds to a mathematical theorem* ...." This insistence that a stochastic process is a mathematical object not to be confused with the empirical process of which it may be a model may seem obvious today but that was not the case in the 1930s.

In order to appreciate Doob's early achievements, one must first have some understanding of the state of probability at the beginning of the 1930s. It was not clear at the time whether probability was part of mathematics or part of physical science. Although many important results had been established, there was no "theory" of probability. Doob, himself, described this state of affairs in [81]. Then in 1933, Kolmogorov proposed the now familiar axiomatic system for probability based on measure theory. Today this is almost universally accepted as the appropriate framework for mathematical probability.

In fact, Kolmogorov showed how to construct a family $\{x_t : t \in T\}$ of random variables with given finite dimensional distributions where $T$ is an arbitrary index set. This Kolmogorov construction was adequate when $T$ is countable; in particular, it provided a mathematical model for a sequence of independent random variables with given distributions. Hence, it was sufficient to formulate and discuss the classical results of probability theory such as the central limit theorem, the strong law of large numbers, etc. However, this is no longer the case when $T$ is uncountable. For example, if $T = \mathbb{R}$, the probabilities that $t \to x_t$ is continuous or that $t \to x_t$ is measurable are not defined in the Kolmogorov construction. But such regularity properties of the trajectories are precisely what are needed in developing a theory of stochastic processes (i.e., family of random variables) indexed by $\mathbb{R}$. Therefore, it is necessary to modify the Kolmogorov construction in order to have the probabilities of such events defined.

Although Wiener had given a rigorous construction of Brownian motion in the 1920s, there was hardly any theory of continuous parameter stochastic processes in the mid 1930s. This is the problem that Doob addressed for the first time in his 1937 paper [12]. He returned to this problem a number of times, especially in [27] and his book [36]. He was the first to rigorously address such questions. This required a serious use of measure theory and I believe it is fair to say that many probabilists in the late 1930s and 1940s (and, perhaps later) were uncomfortable with the type of measure theory required for the analysis of the trajectories of a continuous parameter process. Of course, there were exceptions, especially K. Itô who has written that he was influenced by Doob's work.

Joe Doob was a true pioneer in this area. Elsewhere I have called him a lonely voice in the wilderness during much of this time. (In a letter dated September 27, 1979, he wrote, "Your idea that I was a lonely voice is quite justified.") The appearance of his book *Stochastic Processes* in 1953 was an event of great importance in the theory of probability. It contained along with many other things a culmination of his work on continuous parameter process and a presentation of the state of the art at that epoch of martingale theory. P. A. Meyer has said that it marked the beginning of the modern era in the theory of stochastic processes. The enormous influence of this book cannot be measured by a mere list of citations. It helped shape the spirit and outlook of the generation of probabilists that came of age in the 1950s and 1960s. It is hardly ever cited directly today, having been superseded by a myriad of modern books. It was so successful that it became obsolete.

On a personal note, I first saw, I won't say met, Joe during the summer of 1953. I was in my last year as a graduate student at the University of Michigan and Marcel Brelot was spending a month or two visiting the University of Michigan that summer. At some point during his visit, there was a conference on potential theory and Joe was one of the participants. I must have attended his talk, but I don't recall anything about it. That fall Shu–Teh Moy returned to Ann Arbor to complete her doctorate. She had spent the previous year in Urbana studying with Joe and she brought back some of the page proofs—or perhaps the galleys—of the 1953 book. During the Fall semester of 1953 we organized a seminar to work through the parts of Chapters 1 and 2 dealing with conditioning and separable processes. As there were no copy machines in those days, the proof sheets were just about worn out by the end of the term. The material on separability was especially difficult for me and it was only much later that I understood what he had achieved.

Since his death in June 2004 at the age of 94, numerous testimonials to Doob and his work have appeared. Volume 50 (2006) of the Illinois Journal of Mathematics consists entirely of papers dedicated to his memory and contains a short commentary, by Donald Burkholder, on his life and work. See also [B] and [BP]. I especially recommend [Sn] to the reader.

I am going to attempt to describe Doob's work in two of the areas in which his contributions were seminal: the establishment of a measure theoretic foundation for continuous parameter stochastic processes and probabilistic potential theory. In preparing these remarks I have benefited greatly from and also borrowed freely from two historical notes by P. A. Meyer [M2], [M3] and from [Sh] by M. J. Sharpe.

**2. Continuous parameter stochastic processes.** As mentioned before in his 1933 monograph, Kolmogorov had shown how to construct a family of real valued random variables with prescribed finite dimensional distributions; a construction that is familiar to all probabilists. Let $T$ be an arbitrary index set and let $\Omega^*$ be the set of all real valued functions $\omega^* \colon T \to \mathbb{R}$. Let $\mathcal{B}$ denote the Borel $\sigma$-algebra of $\mathbb{R}$ and let $\mathcal{F}^0$ be the product $\sigma$-algebra in $\Omega^*$, that is, the $\sigma$-algebra generated by the coordinate maps $X_t \colon \Omega^* \to \mathbb{R}$ where $X_t(\omega^*) = \omega^*(t)$. It is convenient to use $X(t, \omega^*)$ and $X_t(\omega^*)$ interchangeably as notation for $\omega^*(t)$. Given a family of finite dimensional distributions subject to the necessary compatibility conditions, Kolmogorov constructed a probability measure $P^*$ on $(\Omega^*, \mathcal{F}^0)$ such that $\{X_t; t \in T\}$ is a family of random variables with the given finite dimensional distributions. Let $\mathcal{F}^*$ denote the $P^*$-completion of $\mathcal{F}^0$.

The intuition is that the finite dimensional distributions associated with an empirical process are the only things that are physically observable, at least theoretically, and so any family of random variables with the given finite dimensional distributions ought to be a reasonable model for the underlying empirical process. This is indeed the case when $T$ is countable. The Kolmogorov space $(\Omega^*, \mathcal{F}^*, P^*)$ provides a satisfactory model for a sequence of random variables with prescribed finite dimensional distributions, in particular, for a sequence of independent identically distributed random variables.

However, if $T$ is uncountable, this is no longer the case. Until further notice, I shall suppose that $T$ is either $\mathbb{R}$ or $\mathbb{R}^+ = [0, \infty[$. A set in $\mathcal{F}^0$ depends on at most a countable number of $t$-values, but in doing analysis one is interested in events such as $\{\omega^* : X(\cdot, \omega^*) \text{ is continuous}\}$ or $\{\omega^* : X(\cdot, \omega^*) \text{ is measurable}\}$ and these sets are not in $\mathcal{F}^0$ nor in $\mathcal{F}^*$. Therefore, their $P^*$ probabilities are not defined. This is the problem that Doob addressed for the first time in his 1937 paper [12]. He began with the Kolmogorov space $(\Omega^*, \mathcal{F}^*, P^*)$ where, as above, $\mathcal{F}^*$ is the $P^*$-completion of $\mathcal{F}^0$. If $\Lambda \subset \Omega^*$, define $\overline{P^*}(\Lambda)$ as the infimum of $P^*(\Lambda')$ for $\Lambda' \in \mathcal{F}^*$ with $\Lambda' \supset \Lambda$. $\overline{P^*}$ is the outer measure induced by $P^*$. If $\Omega \subset \Omega^*$ with $\overline{P^*}(\Omega) = 1$, define a $\sigma$-algebra $\mathcal{F}$ on $\Omega$ to consist of those sets $\Lambda \subset \Omega$ such that $\Lambda = \Omega \cap \Lambda^*$ for some $\Lambda^* \in \mathcal{F}^*$ and set $P(\Lambda) = P^*(\Lambda^*)$. Doob first proved that $P$ is well defined and that $(\Omega, \mathcal{F}, P)$ is a probability space. $\mathcal{F}$ (resp. $P$) is called the *trace* of $\mathcal{F}^*$ (resp. $P^*$) on $\Omega$. Although this is easily checked, I believe that it had not been formalized previously. In any case, he certainly was the first to appreciate its importance. Indeed, he proceeded to *define* a stochastic process as a subset $\Omega$ of $\Omega^*$ with outer measure one. Since $(\Omega^*, \mathcal{F}^*, P^*)$ is complete, so is $(\Omega, \mathcal{F}, P)$. In view of his later insistence on not confusing the mathematical object with the empirical process of which it may be a model, it is interesting that he qualifies this definition with the following footnote: "Strictly speaking, the stochastic process should be defined as the physical system or other entity whose changing is represented by the mathematical formulation of the definition, but it seems wiser to use the term stochastic process both for the mathematical formulation and for the concretization it represents, than to introduce more terminology." This definition suffers from the fact that it involves the rather arbitrary choice of the set $\Omega$. Doob, himself, makes this comment in [27].

In [12] he introduced two concepts. First of all, a process $\Omega$ is *quasi-separable* provided there exists a countable dense set $S \subset T$ such that, for every open interval $I \subset T$,

$$\sup_{t \in I \cap S} X_t(\omega) = \sup_{t \in I} X_t(\omega) \quad \text{and} \quad \inf_{t \in I \cap S} X_t(\omega) = \inf_{t \in I} X_t(\omega)$$

for all $\omega \in \Omega$. Actually his definition is slightly different; the above is an equivalent formulation as is easily seen and is contained in the proof of his Theorem 1.2. One can show that $\Omega$ is quasi-separable if and only if for each $\omega \in \Omega$ the graph of $\omega$ restricted to the countable set $S$ is dense in the graph of $\omega$. Since the graph of any real valued function on $T$ is separable as a subspace of $T \times \mathbb{R}$, quasi-separability is a type of "uniform" separability of the graphs of $\omega \in \Omega$ in that a single countable parameter set $S$ suffices for all $\omega \in \Omega$. Note that if $\Omega$ is quasi-separable, then $\sup_{t \in I} X_t$ as a function on $\Omega$ is $\mathcal{F}$ measurable since it is identical with $\sup_{t \in I \cap S} X_t$ on $\Omega$; the crucial point being that this and many other functions on $\Omega^*$ that are not $\mathcal{F}^*$ measurable become $\mathcal{F}$ measurable when restricted to $\Omega$.

Doob defined a stochastic process, $\Omega$, to be measurable provided that $(t, \omega) \to X(t, \omega)$ is measurable on $T \times \Omega$ relative to the product measure $\lambda \times P$ where $\lambda$ is a Lebesgue measure on $T$. By this he meant measurable with respect to the $\lambda \times P$ completion of $\mathcal{B} \times \mathcal{F}$ or, equivalently, the $\lambda \times P$ completion of $\mathcal{L} \times \mathcal{F}$. Here $\mathcal{B}(\mathcal{L})$ is the $\sigma$-algebra of the Borel (Lebesgue) measurable subsets of $T$. He showed that the set of $\omega^*$ such that $X(\cdot, \omega^*)$ is *not* Lebesgue measurable always has $P^*$ outer measure one for any Kolmogorov space. This and Fubini's theorem imply that $\Omega^*$ itself is never a measurable stochastic process.

In Theorem 2.3 he obtained a necessary and sufficient condition on $P^*$ measure that there exists a measurable process $\Omega$ in $\Omega^*$. Perhaps more interesting is the sufficient condition contained in Theorem 2.4. Namely, if for $\lambda$ almost every $t, X_{t+h} \to X_t$ in $P^*$ measure as $h \to 0$, then there exists a measurable process $\Omega \subset \Omega^*$. He also established a sufficient condition for the existence of a quasi-separable measurable process $\Omega$ in Theorem 2.5. I won't describe this since better results are contained in [16]. Most of the remainder of [12] is devoted to establishing the basic properties of processes with independent increments, here called differential processes. It culminates with the following result which, for simplicity, I shall state in modern terminology. A centered process with independent increments has a version in which almost surely the sample functions have only discontinuities of the first kind and are right continuous except possibly, at the fixed points of discontinuity which form a countable set.

I now turn to the 1940 paper [16]. This paper is best known for its development of the basic properties of martingales including the convergence theorems and the sample function properties of continuous parameter martingales. But Section 2 contains an expanded and improved treatment of much of the material in [12]. Adopting the same definition as in [12], he proved that there *always* exists a quasi-separable process $\Omega$, although in some cases it may be necessary to allow the functions in $\Omega$ to take infinite values. However, the coordinate functions $(X_t)$ have the same finite dimensional distributions on $(\Omega, \mathcal{F}, P)$ as they do on $(\Omega^*, \mathcal{F}^*, P^*)$. Of course this is essentially the basic theorem for the existence of a separable version as formulated in the 1953 book [36]. (I do not know if [36] is the first place where he adopted the terminology separable in place of quasi-separable.) He also showed that if a measurable process exists, then a quasi-separable measurable process $\Omega$ exists which, in fact, possesses additional nice properties that I won't describe here. For later applications the following result is, perhaps, the most important. It is contained in his Theorem 2.8. Suppose that $P^*$ has the property that if $S$ is any countable dense subset of $T$, then almost surely $P^*$, $\lim_{s \in S, s \downarrow t} X_s$ exists for each $t$. Let $D_1(D_2)$ be the subsets of $T$ such that, for any sequence, $(t_n) \subset T$ with $t_n \downarrow t (t_n \uparrow t), X_{t_n} \to X_t$ almost surely $P^*$ provided $t \notin D_1(D_2)$. Then $D_1$ and $D_2$ are countable. Moreover, there exists a quasi-separable process, $\Omega$, such that, for almost all $\omega \in \Omega, X(\cdot, \omega)$ has a limit from the right at all $t$, has at most a countable set of discontinuities (depending on $\omega$) and is right continuous at all $t \notin D_1$.

That same year, 1940, in a joint paper [15] Doob and Ambrose discussed the relationship between a stochastic process defined as a subset $\Omega \subset \Omega^*$ of outer measure one and a stochastic process defined as one usually sees it today, namely, a stochastic process defined as a family $\{Y_t, t \in T\}$ of random variables on some probability space $(W, \mathcal{G}^0, Q)$. They called it a process in the sense of Wiener. Let $Y(t, w)$ also denote $Y_t(w)$. The function $(t, w) \to Y(t, w)$ is a random function in their terminology. They began by describing the, now familiar, way that a random function induces a probability $P^*$ on $(\Omega^*, \mathcal{F}^0)$. Recall that $\mathcal{F}^0$ is the $\sigma$-algebra generated by $(X_t; t \in T)$. Let $\overline{Y}$ be the map from $W$ to $\Omega^*$ defined by $\overline{Y}(w)(\cdot) = Y(\cdot, w)$. $\overline{Y}$ is measurable, that is, $\overline{Y}^{-1}(\mathcal{F}^0) \subset \mathcal{G}^0$. Define $P^*$ to be the image of $Q$ under $\overline{Y}, P^*(\Lambda) = Q[\overline{Y}^{-1}(\Lambda)]$ for $\Lambda \in \mathcal{F}^0$, and $\mathcal{F}^*$ to be the $P^*$ completion of $\mathcal{F}^0$. Then $(\Omega^*, \mathcal{F}^*, P^*)$ is a Kolmogorov space on which the coordinate maps $(X_t; t \in T)$ have the same finite dimensional distributions as $(Y_t, t \in T)$ on $(W, \mathcal{G}^0, Q)$. Set $\Omega = \Omega(Y) = \overline{Y}(W) \subset \Omega^*$. Then as in [12] or [16], $(\Omega, \mathcal{F}, P)$ denotes the trace of $(\Omega^*, \mathcal{F}^*, P^*)$ on $\Omega$. Also let $\mathcal{G}$ be the $\sigma$-algebra on $W$ generated by $(Y_t; t \in T)$ and $\overline{\mathcal{G}}$ its $Q$-completion. Now $\mathcal{G} \subset \mathcal{G}^0$ and the inclusion may be proper. Next they proved that $\Lambda$ is $P$-measurable (i.e., $\Lambda \in \mathcal{F}$) if and only if $\Lambda' = \overline{Y}^{-1}(\Lambda) \in \overline{\mathcal{G}}$ and, hence, $P(\Lambda) = Q(\Lambda')$. They remarked that if $(\Omega^*, \mathcal{F}^*, P^*)$ is a Kolmogorov space, then any stochastic process $\Omega \subset \Omega^*$ is determined by a random function, namely, $X(t, w)$ on $(\Omega, \mathcal{F}, P)$. Therefore, there is an equivalence of the two ways of defining a stochastic process.

Let $Y$ be a random function on $(W, \mathcal{G}^0, Q)$ and, as in the previous paragraph, $\mathcal{G}$ be the $\sigma$-algebra generated by $(Y_t, t \in T)$. They defined $Y$ to be *strongly measurable* provided $(t, w) \to Y(t, w)$ is measurable with respect to the $\lambda \times Q$ completion of $\mathcal{L} \times \mathcal{G}$. They proved that $\Omega(Y) \subset \Omega^*$ is a measurable stochastic process if and only if $Y$ is strongly measurable. Two random functions $Y$ and $Z$ on the same probability space $(W, \mathcal{G}^0, Q)$ are equivalent provided that $Y(t, \cdot) = Z(t, \cdot)$ a.e. $Q$ for each fixed $t \in T$. Then the following result was established. A random function is strongly measurable provided it is equivalent to a random function that is $(t, w)$ measurable, that is, measurable with respect to the $\lambda \times Q$ completion of $\mathcal{L} \times \mathcal{G}^0$. This differs from strongly measurable in that $\mathcal{G}$ is replaced by the possibly larger $\sigma$-algebra $\mathcal{G}^0$. As a corollary, if $Y$ is a random function and $P^*$ is the measure on $\Omega^*$ that it induces, then $\Omega^*$ will contain a measurable process $\Omega$ provided $Y$ is equivalent to a $(t, w)$ measurable random function.

In [27], which is based on an address delivered at the Summer Meeting of the AMS in August 1946, Doob proposed an alternative approach. He also changed his notation somewhat in [27], and I will adapt my notation to conform more closely with that of [27]. Given a family of compatible finite dimensional distributions, Doob now let $(\Omega, \mathcal{F}_0, P_0)$ denote the corresponding Kolmogorov space except that each $\omega \in \Omega$ is now a function from $T$ to $\overline{\mathbb{R}} = [-\infty, \infty]$, that is, $\Omega = \overline{\mathbb{R}}^T$. Of course, $P_0(|X_t| < \infty) = 1$ since the finite dimensional distributions are carried by $\mathbb{R}^d$ for the appropriate $d < \infty$. Note that $\Omega$ equipped with the product topology is a compact Hausdorff space. Let $\Omega_1 \subset \Omega$ have outer measure one. Doob asserted that the method of "relative measure," that is, taking the trace on $\Omega_1$, is inelegant and that it is more elegant to extend the domain of the definition of $P_0$. So let $\mathcal{F}_1$ be the $\sigma$-algebra in $\Omega$ generated by $\mathcal{F}_0$ and the set $\Omega_1$ of outer measure one. It is easy to see that $P_0$ extends uniquely to a probability $P_1$ on $\mathcal{F}_1$. Clearly $(X_t, t \in T)$ has the same finite dimensional distributions when considered on either the trace of $\Omega$ on $\Omega_1$ or on $(\Omega, \mathcal{F}_1, P_1)$, and so this procedure accomplishes exactly what the trace does. This allowed him to more easily compare this technique to a method suggested by Kakutani, namely, to extend $P_0$ on $\mathcal{F}_0$ to $P_2$ defined on the Borel $\sigma$-algebra $\mathcal{F}_2$ of the compact Hausdorff space $\Omega$. Since $\mathcal{F}_0$ is the Baire $\sigma$-algebra of $\Omega$, that is, the $\sigma$-algebra generated by the continuous functions, there is a unique regular extension $P_2$ of $P_0$ to $\mathcal{F}_2$. P. A. Meyer has called $(\Omega, \mathcal{F}_2, P_2)$ the second canonical process in [M1]. Since $\{\sup_{t \in I} X_t > k\}$ is open for any $I \subset T$, this supremum is a random variable on the probability space $(\Omega, \mathcal{F}_2, P_2)$. Note that $\mathcal{F}_1 \subset \mathcal{F}_2$ provided the additional set $\Omega_1 \in \mathcal{F}_2$. But this is still not satisfactory. Doob showed, for example, that there are difficulties when discussing measurable processes. For example, if $\Omega_1$ is a measurable process as defined in [12], $\Omega_1$ may not be in $\mathcal{F}_2$.

In Sections 1 and 2 of Chapter II of his 1953 book [36], Doob gave an expanded treatment of this material. Here he adopted as the definition of a stochastic process what in [15] was called a process in the sense of Wiener. Namely, an arbitrary probability space $(\Omega, \mathcal{F}, P)$ is given, but it is no longer assumed that $\mathcal{F}$ is complete. A random variable $X$ is an extended real valued measurable function on $\Omega$ with the property that $P(|X| < \infty) = 1$, and a stochastic process is a collection $(X_t, t \in T)$ of random variables where $T$ is an arbitrary index set. As before, $T$ will denote either $\mathbb{R}$ or $\mathbb{R}^+ = [0, \infty[$ in what follows and either $X_t(\omega)$ or $X(t,\omega)$ will be used to denote the evaluation of $X_t$ at $\omega \in \Omega$. Parenthetically, it is interesting to note that here he used the name *random variable,* while in the earlier papers [12, 15, 16, 27] he used *chance variable.* On page 307 of [Sn] he explained that this change of nomenclature was the result of a chance (random) experiment: He lost a coin toss. Of course, today random variable is the well established terminology and I, personally, am glad that he lost.

Let $X = (X_t, t \in T)$ be a stochastic process and $\mathcal{A}$ a subfamily of the $\sigma$-algebra, $\mathcal{B}$, of Borel subsets of $T$. Then $X$ is *separable relative* to $\mathcal{A}$ if there is a countable subset $S \subset T$ and $\Lambda \in \mathcal{F}$ with $P(\Lambda) = 0$ such that if $A \in \mathcal{A}$ and $I$ is any open interval, then $\{X_t \in A, t \in I \cap T\}$ and $\{X_s \in A, s \in I \cap S\}$ differ by a subset of $\Lambda$. If $\mathcal{A}$ is the class of all (possibly unbounded) closed intervals, then $X$ is just *separable* without reference to $\mathcal{A}$. If $\Lambda$ is empty, separable is essentially what was called quasi-separable in [12]. In this situation $X$ is separable if and only if the graph of $X(\cdot, \omega)$ is contained in the closure in $T \times \overline{\mathbb{R}}$ of the graph of $X(\cdot, \omega)$ restricted to $S$ for each $\omega \notin \Lambda$. Let $\mathcal{F}(T)$ be the $\sigma$-algebra generated by $(X_t, t \in T)$ and $\overline{\mathcal{F}}(T)$ be the completion of $\mathcal{F}(T)$ relative to $P$. A second process $\widetilde{X} = (\widetilde{X}_t, t \in T)$ on the same probability space $(\Omega, \mathcal{F}, P)$ is a standard modification of $X$ provided that, for each $t \in T, \{X_t \neq \widetilde{X}_t\} \in \overline{\mathcal{F}}(T)$ and $P(X_t \neq \widetilde{X}_t) = 0$. This is stronger than just requiring that $P(X_t \neq \widetilde{X}_t) = 0$ for $t \in T$ since the exceptional set for each $t$ depends only on the process $X$.

THEOREM 1.    *Let $X$ be a stochastic process. Then there exists a standard modification $\widetilde{X}$ of $X$ that is separable relative to the class of all closed sets, and hence, separable.*

This is a stronger result than the corresponding theorem in [16].

The definition of a measurable process in the book is that a process $(X_t, t \in T)$ is measurable provided $X_t(\omega)$ defines a function measurable in the pair of variables $(t, \omega)$. Reading the proof of the following theorem, it is evident what is meant here is that $(t, \omega) \to X_t(\omega)$ is measurable with respect to the $\lambda \times P$ completion of $\mathcal{B} \times \mathcal{F}(T)$. Here is the basic existence theorem for a standard modification given in the book.

THEOREM 2.    *Let $(X_t, t \in T)$ be a stochastic process. If there exists $T_1 \subset T$ with $\lambda(T_1) = 0$ such that $\lim_{s \to t} X_s = X_t$ in probability for each $t \in T \setminus T_1$, then there exists standard modification of $X$ that is measurable and separable relative to the closed sets.*

On page 67 (1953 edition) Doob states a somewhat complicated necessary and sufficient condition for a process $X$ to have a measurable standard modification. Then he goes on to say (emphasis added), "Moreover, it can be shown that there is a separable (relative to the closed sets) measurable standard modification whenever there is a measurable standard modification. *Note that the stated necessary and sufficient condition is a condition on the measures of $\mathcal{F}(T)$ sets, that is, on the finite dimensional distributions.*"

The two theorems quoted above were the culmination of Doob's efforts to develop a framework general enough for a rigorous treatment of continuous parameter processes. In the late sixties and early seventies several new types of process measurability appeared. Predictable, optional or progressive processes were defined and proved to be of great importance. For the most part, these ideas were introduced and developed by P. A. Meyer and his school and became known as the general theory of processes or, sometimes, Strasbourg probability. However, progressive or progressively measurable processes were defined by Chung and Doob in [67] where they also established the basic properties of such processes. Doob returned to this subject in [85] where he related the older notion of separability to these newer (in 1975) concepts.

Let $(\Omega, \mathcal{F}, P)$ be a complete probability space and let $(\mathcal{F}_t, t \geq 0)$ be a right continuous filtration such that $\mathcal{F}_0$ contains all $P$ null sets. Then if $0 \leq s < t < \infty, \mathcal{F}_s \subset \mathcal{F}_t$ and $\mathcal{F}_s = \bigcap_{t>s} \mathcal{F}_t$. A process $X = (X_t, t \geq 0)$ will always have $\mathbb{R}^+$ as parameter set and be adapted to $(\mathcal{F}_t)$. A stopping time is a stopping time relative to the given filtration, and a stopping time is *discrete* provided its range is countable. In [85] Doob proved the following result. Let $X$ be a bounded separable process. Then a sufficient condition that almost surely $X$ has left (right) limits at all $t > 0$ $(t \geq 0)$ is that whenever $(T_n)$ is a uniformly bounded increasing (decreasing) sequence of discrete stopping times $\lim_n E(X(T_n))$ exists. This should be contrasted with the more familiar result that this holds with separable replaced by optional. Of course, every process, $X$, has a separable modification, $\widehat{X}$. Since $P(X_t \neq \widehat{X}_t) = 0, \widehat{X}$ is also adapted to $(\mathcal{F}_t)$. Recall $\mathcal{F}_0$ contains all $P$ null sets.

Doob also defined two new types of separability in [85]. A process $X$ is optionally (predictably) separable provided there exists a sequence $(T_n)$ of bounded (predictable) stopping times such that for each $\omega$ the set $\{T_n(\omega), n \geq 1\}$ contains 0, is dense in $[0, \infty[$ and the graph of $X(\cdot, \omega)$ is in the closure of the graph restricted to this parameter set. He then proved that if $X$ is optional (predictable), it is indistinguishable from an optionally (predictably) separable process $\widehat{X}$. The result described in the previous paragraph remains true if separable is replaced by optionally separable except that one can no longer restrict the $T_n$ to be discrete. There is also a predictable version. If $X$ is a bounded predictably separable process and if $\lim_n E(S(T_n))$ exists whenever $(T_n)$ is a uniformly bounded increasing sequence of predictable stopping times, then almost surely $X$ has left limits.

Today Doob's approach to continuous parameter processes based on separability, for the most part, has been replaced by the newer concepts mentioned above, even though separability is more general in that every process has a separable version. But it was Doob's work during the fifteen years leading up to the publication of his book [36] that set the stage for the major developments in the theory of continuous parameter stochastic processes that were to come in the years immediately following the appearance of [36]. Moreover, the principal probabilists, at least in the West, responsible for this sudden burst of progress in the late 1950s and 1960s belonged to the generation steeped in the rigorous approaches of his book.

**3. Potential theory.** Doob's first paper on potential theory [37] appeared in 1954 just one year after his book [36]. This marked a shift in his research interests, as most of his efforts during the next few years were devoted to what is sometimes called probabilistic potential theory. Although his thesis and first few papers dealt with analytic function theory, he had not published in this area for some time. Yet in [37] he made use of some of the recent (in 1954) deep results of H. Cartan and M. Brelot in potential theory. A major theme underlying much of Doob's work in potential theory is that a superharmonic (harmonic) function composed with Brownian motion yields an almost surely *continuous* supermartingale (martingale) subject to finiteness conditions, and that this provides a powerful tool for investigating properties of such functions, especially their boundary behavior. Here the key word is continuous. That a superharmonic function composed with Brownian motion is a supermartingale results from a rather trite calculation. That this composition is almost surely continuous is a deep result that is crucial for the applications. It is one of the most important results in Doob's early work in potential theory.

Recall that $\Delta = \sum_{j=1}^{d} \frac{\partial^2}{\partial x_j^2}$ denotes the Laplacian in $\mathbb{R}^d$. Let $|x - y|$ denote the Euclidean distance between $x$ and $y$ in $\mathbb{R}^d$. Then $B(x,r) = \{y \in \mathbb{R}^d : |x - y| < r\}$ and $S(x,r) = \{y \in \mathbb{R}^d : |x - y| = r\}$ denote the open ball and sphere with center $x$ and radius $r$. If $G \subset \mathbb{R}^d$ is open, then $u : G \to \overline{\mathbb{R}} = [-\infty, \infty]$ is *superharmonic* on $G$ provided $u > -\infty, u \not\equiv +\infty, u$ is lower semicontinuous (lsc) and for every ball whose closure $\overline{B}(x,r) \subset G$ one has $u(x) \geq \int_{S(x,r)} u(y)\sigma_{x,r}(dy)$, where $\sigma_{x,r}$ is normalized surface measure on $S(x,r)$. Also $u$ is subharmonic provided $-u$ is superharmonic and $u$ is harmonic provided it is both superharmonic and subharmonic. Of course, $u$ is harmonic if and only if it is twice differentiable and $\Delta u = 0$. In fact, a harmonic function is real analytic.

Let $B = (B_t)_{t \geq 0}$ denote Brownian motion in $\mathbb{R}^d$ and let $P^x$ denote its law starting from $x \in \mathbb{R}^d$; $P^x(B_0 = x) = 1$. Let $(\Omega, \mathcal{G})$ be the sample space on which $B$ is defined. We write $B_t(\omega)$ or $B(t,\omega)$ when we want to display the dependence on the sample point $\omega \in \Omega$. It sometimes is also convenient to write $B(t)$ or $B(t, \cdot)$ for $B_t$.

If $f$ is a Borel function on $\mathbb{R}^d$ that is positive (i.e., non-negative) or bounded, define

$$(1) \qquad P_t f(x) = E^x[f(B_t)] = \int g_t(y-x) f(y)\, dy.$$

Here $dy$ denotes $d$-dimensional Lebesgue measure and $g_t(x) = (2\pi t)^{-d/2} \exp(-|x|^2/2t)$ is the Gauss kernel on $\mathbb{R}^d$. Setting $P_0 = I$, one checks that $(P_t)_{t \geq 0}$ is a strongly continuous semigroup of contraction operators on $C_0(\mathbb{R}^d)$ – the continuous functions on $\mathbb{R}^d$ vanishing at infinity.

Some superficial relations between Brownian motion and classical potential theory are immediate. For example, $\frac{1}{2}\Delta$ is the formal generator of the semigroup $(P_t)$, $g_t(x)$ is the fundamental solution of the heat equation $\frac{\partial u}{\partial t} = \frac{1}{2}\Delta u$ and if $d \geq 3$, $\int_0^\infty g_t(x)\, dt = c(d)|x|^{2-d}$ is a multiple of the Newtonian potential kernel. In 1944 Kakutani [K] took the first nontrivial step in a development that would come to recognize that as M. J. Sharpe has written in [Sh] "Brownian motion and classical (Newtonian) potential theory are just different mechanisms for studying the same problems." He identified the Brownian exit distribution from a domain with the corresponding harmonic measure. To explain this, suppose that $A$ is a Borel subset of $\mathbb{R}^d$, then the hitting time, $T(A)$, of $A$ is the infimum of the *strictly* positive $t > 0$ such that $B_t \in A$ and $T(A) = \infty$ if there are no such $t$. Here the fact that the infimum excludes $t = 0$ is crucial. The exit time, $\tau(A)$, from $A$ is the hitting time of $A^c = \mathbb{R}^d \backslash A$. Now let $D$ be a connected open set with compact closure and boundary $\partial D$. If $x \in D$, $H_D(x, dy) = P^x[B_{\tau(D)} \in dy]$ defines a probability measure on $\partial D$. Here $B_{\tau(D)}$ is the position of the Brownian trajectory when it first exits $D$. If $\partial D$ is smooth and $f$ is continuous on $\partial D$, Kakutani showed that

$$(2) \qquad u(x) = \begin{cases} \int H_D(x, dy) f(y), & x \in D, \\ f(x), & x \in \partial D \end{cases}$$

solves the Dirichlet problem for $D$ with boundary data $f$ in the sense that $u$ is harmonic in $D$, continuous on $\overline{D}$ and $u = f$ on $\partial D$. In other words, $H_D(x, \cdot)$ is a harmonic measure for $D$. Kakutani did not give complete proofs of either the harmonicity of $u$ in $D$ or the continuity of $u$ on $\overline{D}$. The proof that $u$ is harmonic in $D$ required a special case of the strong Markov property (*SMP*). Kakutani clearly stated the result that was needed, but omitted the proof since it "requires a complicated argument."

The next major advance was Doob's paper [37]. This paper might best be described as an application of martingale theory to the study of harmonic and subharmonic functions. After pointing out that there is a formal parallelism between subharmonic functions and submartingales, he wrote,

"We are more interested in a study which brings the two topics together, and thereby obtain new results in the theory of subharmonic and harmonic functions." This parallelism is explicitly spelled out in his 1984 book *Classical Potential Theory and Its Probabilistic Counterpart* [91]. See also the expository paper [89]. Doob considered subharmonic functions but I shall state his results for superharmonic functions. In the body of [37] he supposed that $d = 2$, but the results carry over to $d > 2$ with only minor changes as he remarked in a concluding section, and I shall suppose $d \geq 2$ in what follows. (The results when appropriately interpreted also hold for $d = 1$.)

In [37] Doob proved a special case of Kakutani's result "since Kakutani has not yet published the details of his proof." Namely, if $0 < r < |x| < R < \infty$, then the probability $p(x)$ that Brownian motion starting from $x$ hits $S_R = \{y : |y| = R\}$ before it hits $S_r = \{y : |y| = r\}$ is the unique harmonic function in the annulus $\{x : r < |x| < R\}$ with boundary values one on $S_R$ and zero on $S_r$. Of course, this harmonic function has a simple explicit expression depending on the dimension. In proving this, he used $(SMP)$ without explicitly mentioning it. But earlier in the paper he had clearly stated $(SMP)$ and treated it with the following comment: "Although this fact gives considerable insight into some of the results of this paper, it will not be proved, because it is intuitively obvious, and because it is a special case of a much more general theorem which will be proved elsewhere." It is interesting that this evaluation of $p(x)$ is the only place in the paper that $SMP$ is used. The first proofs of $SMP$ for Brownian motion (and more general processes) appeared almost simultaneously in 1956 in the United States [H] and in the Soviet Union [DY]. However, Doob had stated and proved $SMP$ for a class of continuous parameter Markov chains in 1945 [26]. Just before its statement he wrote, "The following theorem, which is a *special case of a very general theorem on Markoff processes...*" (emphasis added).

Doob went on to show that if $u$ is superharmonic (harmonic) on $\mathbb{R}^d$, then, subject to an integrability assumption, the composition $u(B_t)$ is a supermartingale (martingale). This is an easy consequence of the rotational symmetry of the Gauss kernel. But to turn it into a powerful tool, one needs to establish some regularity properties of $t \to u(B_t)$. Although a superharmonic function $u$ may be very discontinuous and take the value $+\infty$, he proved that $t \to u(B_t)$ is finite and continuous on $]0, \infty[$ almost surely $P^x$ for all $x$, and if $u(x) < \infty$, then zero may be included. To my mind, this is perhaps the most important probabilistic result in the paper, and it certainly is the key technical fact that is needed for some of the deeper results of the paper. Doob based his proof on a previous result of Cartan that states, roughly speaking, given $\varepsilon > 0$, there exists an open set of capacity less than $\varepsilon$ such that $u$ restricted to its complement is continuous.

I have slurred over a very important point. What Doob needed and proved is that if $G$ is an open subset of $\mathbb{R}^d$, $u$ is superharmonic (harmonic) on $G$ and $D$ is open with $\overline{D} \subset G$, then $X_t = u(B_{t \wedge \tau(D)})$ is an almost surely continuous supermartingale (martingale), subject to a finiteness condition, with respect to $P^x$ for $x \in D$. The fact that $X = (X_t; t \geq 0)$ is a supermartingale (martingale) is no longer a simple calculation. Doob proved this by reducing the general case to that in which $G$ is all of $\mathbb{R}^d$. It was in accomplishing this reduction that he needed the special case of Kakutani's result mentioned above.

Doob used these results to obtain an intuitive probabilistic characterization of "thinness," a concept that had been introduced by Brelot. Namely, a set $A$ is thin at $x$ if and only if almost surely, the Brownian path starting from $x$ avoids $A$ on some interval $]0, \varepsilon[$ where $\varepsilon > 0$ depends on $\omega$. Actually there is a zero–one law and one has the criteria that $A$ is thin at $x$ if $P^x(T(A) > 0) = 1$ and $A$ is not thin at $x$ if $P^x(T(A) = 0) = 1$. This holds for any Borel or even analytic set $A \subset \mathbb{R}^d$, although Doob only considered the case in which $A$ is an $F_\sigma$ set. The general result uses Choquet's theorem that analytic sets are capacitable, which only appeared a year later in 1955.

Let $G$ be an open subset of $\mathbb{R}^d$ and, for simplicity, suppose that $G$ has compact closure. There is a class of measurable functions (not necessarily finite) on $\partial G$ such that for each $f$ in the class the Perron–Wiener–Brelot (PWB) method assigns a unique harmonic function $u_f$ on $G$. Such functions are called resolutive for $G$. In particular, continuous functions are resolutive. A point $y \in \partial G$ is regular for the Dirichlet problem provided $\lim_{x \in G, x \to y} u_f(x) = f(y)$ for all continuous $f$ on $\partial G$. In 1940 Brelot had proved that $y \in \partial G$ is regular for the Dirichlet problem on $G$ if and only if $G^c$ is not thin at $y$. Therefore, the probabilistic description of thinness gave the following intuitive condition that $y$ is regular if and only if $P^y(\tau(G) = 0) = 1$. Recall that $\tau(G) = \inf\{t > 0 : B_t \notin G\}$. Now for each $x \in G$, $f \to u_f(x)$ is a positive linear functional on $C(\partial G)$. Therefore, the Riesz representation theorem implies that for each $x \in G$ there exists a probability measure $\mu_G(x, \cdot)$ on $\partial G$ such that $u_f(x) = \int f(y) \mu_G(x, dy)$ for $f \in C(\partial G)$. Thus, $\mu_G(x, \cdot)$ is harmonic measure for $G$. Brelot had established that a necessary and sufficient condition that $f$ be resolutive for $G$ is that $f$ be $\mu_G(x, \cdot)$ integrable for at least one $x$ in each connected component of $G$. Then $f$ is $\mu_G(x, \cdot)$ integrable for all $x \in G$ and $u_f(x) = \int f(y) \mu_G(x, dy)$.

This raises the question: In what sense is a resolutive $f$ the boundary function of $u_f$? Doob gave an elegant probabilistic answer to this question:

THEOREM.  *Let $f$ be resolutive for $G$ and let $u = u_f$. Then for each $x \in G$:*

(i) $\lim_{t \uparrow \tau(G)} u(B_t) = f(B_{\tau(G)})$, *a.s. $P^x$,*

(ii) $u(x) = E^x[f(B_{\tau(G)})]$.

Note the second assertion immediately implies that $\mu_G(x, dy) = P^x(B_{\tau(G)} \in dy)$ and this proves and generalizes Kakutani's identification of harmonic measure. Doob made the following comment about this result: "According to this theorem, the Brownian motion solves the Dirichlet problem whenever this problem has a solution, and gives the simplest possible explicit relation between solution and boundary function."

It is very interesting that in proving (ii) Doob managed to finesse the direct use of the strong Markov property. He began by defining $u$ to be the solution of the Dirichlet problem given by the PWB method and then proved that $u$ so defined satisfied (ii) and also (i). However, it was here that he used the fact that $X_t = u(B_{t \wedge \tau(D)})$ is a martingale for $D$ open with $\overline{D} \subset G$ which depended on the identification of the probability of hitting $S_R$ before $S_r$ as described above. Thus, the proof of the theorem did depend on $SMP$ albeit indirectly.

The following year (1955), Doob turned his attention to the heat equation in [42]. Here the situation was quite different. In [37], Doob made use of the extensive existing results for the potential theory of the Laplace equation, but the existing potential theory for the heat equation was, at best, rudimentary. Doob's key insight was to treat as a Dirichlet problem what heretofore had been treated as a Cauchy problem. But from a probabilistic viewpoint there is no essential difference between the elliptic and parabolic cases. Doob began by defining parabolic and superparabolic functions and establishing their basic properties. I won't go into a detailed discussion. Suffice it to say that they are analogous to harmonic and superharmonic functions. Let $X = (X_t)$ be a space time Brownian motion; that is, if $X_0 = (x, s) \in \mathbb{R}^d \times \mathbb{R}$, then $X_t = (B_t, s - t) \in \mathbb{R}^d \times \mathbb{R}$ where $B$ is a Brownian motion on $\mathbb{R}^d$ with law $P^x$. Doob proved that the composition $u(X_t)$ of a superparabolic function $u$ and the process $(X_t)$ is almost surely right continuous in $t$. Here there was no previous result on which to base a proof such as Cartan's that he had used in [37]. Rather he devised an argument which is essentially the method Hunt adopted several years later to prove the analogous result for a general class of Markov processes. The argument depended on the strong Markov property. But this time there is no mention of it. He went on to discuss the Dirichlet problem for the heat equation. He proved results for the heat equation that were completely analogous to those that he had proved the previous year for Laplace's equation in [37].

These two papers were seminal in that they contained new methods and techniques and opened new areas for research that would be pursued for many years by probabilists and by Doob himself. The key idea underlying most of the main results of [37] and [42] is that a wide class of functions naturally associated with the process in question (superharmonic, super-parabolic) when composed with the process yield well behaved (continuous or right continuous) supermartingales. The same idea proved to be fundamental in much of the future development of the subject. In my opinion [42] was never properly appreciated; most likely because two years later Hunt presented a theory that contained many of the results of [42] as special cases.

Doob's next major innovation was the introduction of the $h$-transform of Brownian motion in [45]. Brelot had developed a "relative potential theory" in which, roughly speaking, an arbitrary strictly positive harmonic function $h$ plays the role of 1 in the classical theory. In [45] Doob showed that this $h$-transform of Brownian motion (described in the next paragraph) corresponded to Brelot's relative theory just as ordinary Brownian motion corresponded to the classical theory. This transform when applied to more general processes became one of the central techniques in many later developments in probabilistic potential theory.

If $d \geq 3$, let $R$ be an arbitrary nonempty open subset of $\mathbb{R}^d$, or if $d = 2$, suppose that the complement of $R$ has positive capacity. In his 1984 book [91], Doob called such sets greenian. In [45] Doob worked more generally in that $R$ was an arbitrary Green space, but, for simplicity, I shall suppose that $R$ is a greenian subset of $\mathbb{R}^d, d \geq 2$. I'll say a few words about Green spaces later. Let $B = (B_t)$ be Brownian motion on $\mathbb{R}^d$ and let $\tau = \tau(R)$ be the exit time from $R$. Then $X_t = B_t$ if $t < \tau$ is Brownian motion on $R$ and it has "lifetime" $\tau$ which may be finite with positive probability. The previous year Hunt [H] had shown that for any open subset $R \subset \mathbb{R}^d, d \geq 1, X$ has a strictly positive transition density $q(t, x, y)$ with respect to the Lebesgue measure on $R$, finite and continuous for $t > 0, x, y \in R$ which is symmetric in $x$ and $y$. Let $h > 0$ be superharmonic on $R$ and define $q^h(t, x, y) = h(x)^{-1} q(t, x, y) h(y)$ for $t > 0, x, y \in R$ provided not both $h(x)$ and $h(y)$ are infinite and set it equal to 1 if $h(x) = h(y) = \infty$. Then one verifies that $q^h$ is a transition density relative to the Lebesgue measure and there exists a Markov process $X^h$ on $R$ with this transition function. $X^h$ is called the $h$-transform of $X$ and its paths $t \to X_t^h$ are called $h$-paths. The process $X^h$ has lifetime $\tau^h$ which may be finite. Doob began by establishing the basic properties of $X^h$. In particular, almost surely $t \to X_t^h$ is continuous on $[0, \tau^h[$ and if $u$ is superharmonic on $R$, then $u \circ X_t^h$ is almost surely continuous on $[0, \tau^h[$. The proof of this last property is the same as that given for $B$ in [37]. Moreover, $X^h$ has the strong Markov property. Doob called $X^h$ a conditional Brownian motion. I'll come to the reason for this shortly. I shall write $P^{x|h}$ for the law of $X^h$ starting from $x \in R$.

I shall discuss briefly some of the main results of [45] which hopefully will make it clear where the title comes from. A positive superharmonic function $h$ is minimal provided the only positive superharmonic functions that it dominates are multiples of it. Doob's main purpose in this paper was to investigate the "boundary behavior" of ratios $u/h$ where $u$ and $h$ are strictly positive superharmonic functions on $R$. He began by showing $(u/h)(X_t^h)$ has a finite limit almost surely $P^{x|h}$ as $t \uparrow \tau^h$ provided $h(x) < \infty$. If $h$ is minimal harmonic, then this limit is $\inf_{x \in R} u(x)/h(x)$ on almost every such path; in particular, the limit is constant almost surely and does not depend on $x$ in this case.

To investigate this limit behavior further, following Brelot, he introduced the Martin boundary $R'$ of $R$; $R \cup R'$ is compact and $R$ is a dense open subset of $R \cup R'$. Let $g(x, y)$ denote the Green function of $R$. The Martin kernel is defined as follows. Fix $x_0 \in R$. Then

(3)
$$K_{x_0}(x, y) = \begin{cases} g(x, y)/g(x, x_0), & \text{if } x, y \in R, \\ \lim_{z \to x, z \in R} K_{x_0}(z, y), & \text{if } x \in R', \, y \in R, \end{cases}$$

the limit existing by the definition of the Martin boundary, and one has the Martin representation: if $u$ is a positive superharmonic function on $R$, then there exists a finite measure $\mu^u$ on $R \cup R'$ such that

(4)
$$u(y) = \int K_{x_0}(x, y) \mu^u(dx).$$

Restricting $\mu^u$ to $R$ (respectively $R'$) in (4) yields the Riesz decomposition of $u$ into a potential $p(y) = \int_R K_{x_0}(x, y) \mu^u(dx)$ and the largest harmonic minorant of $u$, $v(y) = \int_{R'} K_{x_0}(x, y) \mu^u(dx)$.

Doob showed that if $h > 0$ is harmonic on $R$, then almost every $h$-path from a point of $R$ converges to a point of $R'$; that is, $X_t^h$ converges to a point of $R'$ as $t \uparrow \tau^h$ almost surely $P^{x|h}$ for each $x \in R$. In particular, if $h$ is minimal, almost all $h$-paths converge to the same point of $R'$ called the pole of $h$. Thus, the Martin boundary serves as an "exit boundary" for $X$ on $R$. This idea that the Martin boundary is an exit boundary was extended to more general processes and became an important concept in probabilistic potential theory. More generally, if $h$ is a strictly positive minimal superharmonic function, then almost every $h$-path from a point of $R$ converges to a unique point of $R \cup R'$, the pole of $h$. For example, if $h = g(x, \cdot)$ where $g$ is the Green function for $R$, then $x$ is the pole of $h$, or if $h > 0$ is a minimal harmonic function, its pole is in $R'$ and is called a minimal point of $R'$. In this sense $X^h$ may be thought as $X$ "conditioned" to converge to the pole of $h$ as $t \uparrow \tau^h$ whenever $h$ is a strictly positive minimal superharmonic function on $R$.

Doob then investigated the Dirichlet problem for $h$-harmonic functions on $R$ when $h > 0$. He established results in this context that were analogous to those he had obtained in [37] for the ordinary Dirichlet problem. For this it was necessary that $R$ be equipped with its Martin boundary $R'$ and not its Euclidean boundary. Since a function $v$ on $R$ is $h$-harmonic ($h$-superharmonic) provided $u = hv$ is harmonic (superharmonic), the results for the possible limiting values of $v$ at the boundary $R'$ gave results about the ratio $u/h$. This led him to a far-reaching generalization of the classical Fatou theorem for harmonic or analytic functions on the disk. To describe this, one needs the fine topology. This is the topology on $\mathbb{R}^d$ generated by the class of superharmonic functions and had been introduced by H. Cartan. The fine topology is strictly finer (i.e., strictly weaker) than the Euclidean topology if $d \geq 2$ and is the Euclidean topology if $d = 1$. Somewhat later Brelot observed that $x$ is a fine limit point of $A$ provided $A$ is not thin at $x$. Consequently, in light of Doob's characterization of thinness, a set $A$ is finely open provided the Brownian path from a point in $A$ remains in $A$ for an initial interval $[0, \eta[$ where $\eta > 0$ almost surely. Intuitively a set is finely open provided it "looks open" to a traveler along the Brownian path. L. Naim extended the concept of thinness to points of $R'$ and thereby extended the fine topology to $R \cup R'$. Doob was able to show that fine limits were the same as limits along $h$-paths in the following sense. If $x \in R$, let $h = g(x, \cdot)$ or if $x$ is a minimal point of $R'$, let $h > 0$ be a minimal harmonic function with pole at $x$. Then a function $u$ on $R$ has fine limit $b$ at $x$ if and only if $u$ has limit $b$ along almost all $h$-paths from a point of $R$ to $x$.

I can now state the generalization of the Fatou theorem alluded to above. If $u$ and $h$ are positive superharmonic functions on $R$ and $h > 0$, then $u/h$ has a finite fine limit at almost all points of $R \cup R'$ relative to $\mu^h$—the measure in the Martin representation (4) of $h$. If $h$ is harmonic, $\mu^h$ is carried by the minimal points of $R'$. Two years later in [53], Doob presented a nonprobabilistic proof of this relative Fatou theorem.

As mentioned before, Doob actually supposed the $R$ was a Green space in [45] rather than a greenian set in Euclidean space. Green spaces had been introduced in potential theory by Brelot and Choquet. Roughly speaking, a Green space $R$ is a locally compact, connected Hausdorff space with a countable base that is "locally Euclidean" of a fixed dimension $d \geq 2$ and on which a strictly positive nonconstant superharmonic function exists. A connected greenian subset of $\mathbb{R}^d, d \geq 2$ is a Green space. Doob established the existence and properties of Brownian motion on a Green space in [46]. The construction involved piecing together processes defined locally on coordinate patches to obtain a globally defined process. This technique became an important method for constructing processes. However, from a contemporary point of view, his argument was sketchy and the treatment of certain technical points was rather minimal, but the fundamental idea was clear. This paper was probably the first in which the consideration of processes with finite lifetimes was truly important since the locally defined processes on coordinate patches appeared naturally with finite lifetimes. This also was probably the first place where the device of adjoining an additional point to the state space was used in order to reduce the finite lifetime case to the situation in which the lifetime is infinite.

All of the results that I have discussed (and much more), with the exception of Green spaces, received an expanded treatment in his 1984 book [91]. This book of 800 pages is a veritable cornucopia of information about potential theory, including parabolic potential theory. In contrast to the 1953 book which appeared at the beginning of an era of very rapid development in the theory of stochastic process, this book appeared toward, if not at, the end of an era in which potential theory and especially probabilistic potential had undergone rapid advancement and were major areas of mathematical research. Although these areas are still active, in my opinion, they are now more peripheral to the main body of mathematical research.

**Acknowledgments.** I would like to thank D. Aldous and H. Heyer for carefully reading a preliminary version and for making a number of valuable suggestions.

## REFERENCES

[B] BINGHAM, N. H. (2005). Doob: A half-century on. *J. Appl. Probab.* **42** 257–266. MR2144908

[BP] BURKHOLDER, D. and PROTTER, P. (2005). Joseph Leo Doob, 1910–2004. *Stochastic Process. Appl.* **115** 1061–1072. MR2147241

[DY] DYNKIN, E. and JUSHKEVICH, A. (1956). Strong Markov processes. *Teor. Veroyatnost. i Primenen.* **1** 149–155. MR0088103

[H] HUNT, G. A. (1956). Some theorems concerning Brownian motion. *Trans. Amer. Math. Soc.* **81** 294–319. MR0079377

[K] KAKUTANI, S. (1944). Two-dimensional Brownian motion and harmonic functions. *Proc. Imp. Acad. Tokyo* **20** 706–714. MR0014647

[M1] MEYER, P.-A. (1966). *Probability and Potentials.* Blaisdell, Boston. MR0205288

[M2] MEYER, P.-A. (1968). La théorie générale des processus de Markov à temps continu. Unpublished manuscript.

[M3] MEYER, P.-A. (2000). Les processus stochastiques de 1950 à nos jours. In *Development of Mathematics 1950–2000* ( J. P. PIER, ED.) 813–848. Birkhäuser, Basel. MR1796860

[Sh] SHARPE, M. J. (1986). S. Kakutani's work on Brownian motion. In *Contemporary Mathematicians—Selected Works of S. Kakutani* (R. KALLMAN, ED.) **2** 397–401. Birkhäuser, Basel.

[Sn] SNELL, J. L. (1997). A conversation with Joe Doob. *Statist. Sci.* **12** 301–311. MR1619190

DEPARTMENT OF MATHEMATICS
UNIVERSITY OF CALIFORNIA, SAN DIEGO
9500 GILMAN DRIVE
LA JOLLA, CALIFORNIA, 92093-0112
USA
E-MAIL: rgetoor@ucsd.edu